\documentclass[a4paper,12pt,babel]{article}
\usepackage[russian]{babel}
\usepackage[active]{srcltx}

\usepackage{amssymb}

\newenvironment{proof}
{\par\noindent{\bf Доказательство.}} 
{\hfill$\scriptstyle\blacksquare$}   

\newcounter{eq}
\renewcommand\theeq{\thesection.\arabic{eq}}

\newtheorem{defn}{Определение}
\newtheorem{theo}{Теорема}

\newtheorem{rem}[theo]{Замечание}

\newtheorem{example}{Пример}

\newcommand{\rf}[1]{(\ref{#1})}

\newcommand{\bibl}[1]{\bibitem{#1}}

\newcommand{\la}{\lambda}
\newcommand{\ep}{\varepsilon}
\begin{document}
УДК 517.518.36
\vskip10pt
\begin{center}
\bf{\Large{ Одна теорема единственности аддитивных }}
\end{center}
\begin{center}
\bf{\Large{ функций и  ее приложения к некоторым  }}
\end{center}

\begin{center}
\bf{\Large{   ортогональным рядам }}
\end{center}

\begin{center}
 К.А. Керян
\end{center}

%
\abstract{В работе исследуется вопрос восстановления определенной на $\mathcal P$-ичных параллелепипедах аддитивной функции по ее производной по $\mathcal P$-ичным параллелепипедам. Полученная теорема применена к вопросам единственности кратных рядов по  системам Хаара и Прайса.}
\vskip20pt


Из теории тригонометрических рядов хорошо известно, что из сходимости тригонометрического ряда почти всюду к нулю не следует, что все коэффициенты этого ряда равны нулю (см. \cite{bari}, стр. 804). Такая же ситуация имеет место для других классических систем, например для рядов по системам Хаара, Уолша, Франклина .

В работах \cite{Al}, \cite{GG89} впервые были рассмотрены вопросы единственности для почти всюду сходящихся или суммирующихся тригонометрических рядов. Понятно, что при этом были наложены дополнительные условия на ряд. Здесь мы их формулировать не будем, поскольку результаты настоящей работы касаются рядов по системам Хаара, Уолша, а также их обобщений. Отметим только, что аналогичные вопросы для кратных тригонометрических рядов рассмотрены в работах \cite{GG93},\cite{G94}.

Через $\mu(A)$ обозначается Лебеговая мера множества $A$. Напомним определение $A$-интеграла.
\begin{defn}Функция $f:[0,1]\rightarrow \mathbb{R}$ называется $A$-интегрируемой, если
$$\lim_{\lambda\rightarrow\infty}\lambda\mu\{x\in[0,1];|f(x)|> \lambda\}=0,  $$
 и существует предел
$$\lim_{\lambda\rightarrow\infty}\int_0^1[f(x)]_\lambda dx=:(A)\int_0^1f(x)dx,$$
где
$$[f(x)]_\lambda =
\left\{\begin{array}{ll}
f(x),  \ \hbox{при}\ |f(x)|\leq \lambda\\
0, \ \ \quad \hbox{при}\ |f(x)|> \lambda.
\end{array}
\right.$$
\end{defn}

Для рядов по системе Хаара Г.Г. Геворкяном \cite{GG95}, в частности, доказана следующая

\begin{theo}\label{thGG}
Пусть ряд по системе Хаара
\begin{equation}\label{1}
\sum_{n=1}^\infty a_n\chi_n(x)
\end{equation}
почти всюду суммируется к $f(x)$ и
\begin{equation}\label{2}
\lim_{\lambda\rightarrow\infty}\lambda\mu\{x\in[0,1];S^*(x)>\lambda\}=0,
\end{equation}
где $S^*(x)$ мажоранта частных сумм ряда \rf{1}. Тогда ряд \rf{1} является  рядом Фурье-Хаара в смысле $A$-интегрирования, т.е.
$$a_n=(A)\int_0^1f(x)\chi_n(x)dx.$$
\end{theo}

Понятно, что из теоремы \ref{thGG} следуют:\\
1) если ряд \rf{1}  п.в. сходится к нулю и выполняется \rf{2}, то все коэффициенты $a_n$ равны нулю;\\
2) если ряд \rf{1}  п.в. сходится к интегрируемой по Лебегу функции $f$ и выполняется \rf{2}, то \rf{1} является рядом Фурье-Хаара функции $f$.\\

Теорема \ref{thGG} В. Костиным \cite{K03} была распространена на ряды по обобщенной системе Хаара  и системе Прайса (см. определения  ниже). Допустим система Прайса или обобщенная система Хаара порождена последовательностью $P=\{p_n\}_{n=1}^\infty, $ где $\sup p_n<\infty.$ Положим $m_n=\prod_{i=1}^np_i$ и $S^*(x)=\sup_{n\in\mathbb{N}}|S_{m_n}(x)|.$ Для таких систем В. Костин \cite{K03} доказал следующую теорему.

\begin{theo}\label{thK03}
Пусть последовательность частных сумм $S_{m_n}(x)$ ряда по системе Прайса или по обобщенной системе Хаара почти всюду сходится к функции $f(x)$ и существует некоторая последовательность  $\{\lambda_n\}_{n=1}^\infty$ для которой
$$\lim_{k\rightarrow\infty}\lambda_k\mu\{x\in[0,1];S^*(x)>\lambda_k\}=0.$$
Тогда $n$-ый коэффициент данного ряда восстанавливается   по формуле
$$a_n=\lim_{m\rightarrow\infty}\int_0^1\left[f(x)\overline{\chi_n(x)}\right]_{\lambda^n_m}dx,$$
где $\lambda_m^n=\lambda_m\|\chi_n\|_\infty,$ в случае ряда по обобщенной системе Хаара $\chi_n(x)$ и
$$b_n=\lim_{m\rightarrow\infty}\int_0^1\left[f(x)\overline{\psi_n(x)}\right]_{\lambda_m}dx,$$
в случае ряда по системе Прайса $\psi_n(x)$.
\end{theo}

В настоящей работе мы усиливаем  эту теорему, рассматривая более широкий $A\mathcal H$-интеграл, введеный К. Йонедой \cite{Y I}.\\\\

Дадим некоторые  необходимые определения и введем обозначения. Для комлекснозначной функции
$f(x)$ определим срезку неотрицательной функцией $\lambda(x)$ следующим образом
$$[f(x)]_{\lambda(x)}=\left\{
\begin{array}{cc}
f(x),& |f(x)|\leq \lambda(x)\\
0,&|f(x)|>\lambda(x).
\end{array}
\right.$$
\begin{defn}(\cite{Y I},\cite{Y II})
Пусть $\mathcal H=\{h_n(x)\}_{n=1}^\infty$ последовательность действительнозначных интегрируемых по Лебегу на $[0,1]$ функций такая, что п.в. выполняется
$$
0\leq h_1(x)\leq h_2(x)\leq\dots\leq h_n(x)\leq\dots, \ \lim_{n\rightarrow\infty}h_n(x)=\infty.
$$
Тогда $f(x):[0,1]\rightarrow C$ называется $A\mathcal H$интегрируемой на $[0,1],$ если
$$\lim_{n\rightarrow\infty}\int_{\{x\in[0,1];|f(x)|\geq\alpha h_n(x)\}}h_n(x)dx=0, \ \textrm{для всех}\ \alpha>0$$
и существует предел
$$\lim_{n\rightarrow\infty}\int_0^1[f(x)]_{h_n(x)}dx=:(A\mathcal H)\int_0^1f(x)dx,$$
который и называется  $A\mathcal H$-интегралом функции $f(x).$
\end{defn}

Если $h_n(x)\equiv n,$ то соответственно получится обычный $A$-интеграл, а если же $h_n(x)\equiv \lambda_n,$ где $\Lambda:=\{\lambda_n\}_{n=1}^\infty$ возрастающая к бесконечности последовательность, тогда получается обобщение $A$-интеграла, которое обозначим через $A\Lambda$-интеграл.\\

Аналогично можно обобщить определения $A\mathcal H$ интеграла для комлекснозначной функции $f(x),$
определенной на параллелепипеде $I\subset R^n.$
\begin{rem}Если для функции $f:I\rightarrow C$ существует такая последовательность возрастающих функций $\mathcal H=\{h_n(x)\}_{n=1}^{\infty},$ что
\begin{equation}\label{H-G}
\lim_{n\rightarrow\infty}\int_{\{x\in I;|f(x)|\geq h_n(x)\}}h_n(x)dx=0,
\end{equation}
и существует предел
$$\lim_{n\rightarrow\infty}\int_I[f(x)]_{h_n(x)}dx,$$ тогда существует последовательность функций $\mathcal G$ такая, что функция $f(x)$ $A\mathcal G$ интегрируема и
\begin{equation}\label{H-G2}
\lim_{n\rightarrow\infty}\int_I[f(x)]_{h_n(x)}dx=(A\mathcal G)\int_{I}f(x)dx.
\end{equation}
\end{rem}

Действительно, из условия \rf{H-G} следует, что существует  неубывающая последовательность $\{\alpha_n\}_{n=1}^\infty,$ сходящаяся к $+\infty$ такая, что
\begin{equation}\label{H-G-al}
\lim_{n\rightarrow\infty}\alpha_n\int_{\{x\in I;|f(x)|\geq h_n(x)\}}h_n(x)dx=0.
\end{equation}
Обозначив $g_n(x)=\alpha_n h_n(x),$ получим, что для любого $\alpha>0,$ $\alpha g_n(x)>h_n(x),$ при достаточно больших $n$, и следовательно
$$
0\leq \lim_{n\rightarrow\infty}\int_{\{x\in I;|f(x)|\geq \alpha g_n(x)\}}g_n(x)dx\leq
\lim_{n\rightarrow\infty}\alpha_n\int_{\{x\in I;|f(x)|\geq h_n(x)\}}h_n(x)dx=0.$$
Кроме того, из \rf{H-G-al} имеем
$$\left|\int_I[f(x)]_{h_n(x)}dx-\int_I[f(x)]_{g_n(x)}dx\right|\leq\int_{\{x\in I;|f(x)|\geq h_n(x)\}}g_n(x)dx\rightarrow0,$$
следовательно выполняется \rf{H-G2}.\\

Пусть $P^j=\{p^j_i\}_{i=1}^\infty, 1\leq j \leq d$ последовательности натуральных чисел, отличных от 1 и
$\mathcal P=\{P^j\}_{j=1}^d$. Положим
$m^j_0=1$ и $m^j_k=\prod_{i=1}^kp^j_i.$
Через $\Lambda^d_k$ обозначим множество всех $\mathcal P$-ичных параллелепипедов ранга $k$, т.е.

$$\Lambda^d_k=\left\{I;I=\left[\frac{n_1}{m^1_k},\frac{n_1+1}{m^1_k}\right]\times\ldots\times
\left[\frac{n_d}{m^d_k},\frac{n_d+1}{m^d_k}\right], n_1,\ldots,n_d\in Z\right\},$$
а также будем писать, что $r(I)=k, $ если $I\in\Lambda^d_k.$

Пусть $\Lambda^d=\cup_{k\in Z_+}\Lambda^d_k,$ где $Z_+=\{0\}\cup N.$ Комплекснозначную функцию $\Psi,$ определенную на множестве $\Lambda^d,$ назовем аддитивной функцией
на $\mathcal P$-ичных параллелепипедах, если для всех $I, I_1,\ldots, I_n\in\Lambda^d,$ для которых $I=\cup_{i=1}^nI_i$ и $ \hbox{int}\,(I_i)\cap \hbox{int}\,( I_j)=\emptyset, i\neq j,$ имеет место
$$\Psi(I)=\sum_{i=1}^n\Psi(I_i).$$

Точку $x\in R$ назовем иррациональной относительно последовательности $\{p_i\}_{i=1}^\infty$, если $x\cdot \prod_{i=1}^np_i\not\in Z,$ для всех $n\in N.$
Точку $x=(x_1,\ldots, x_d)\in R^d$ назовем $\mathcal P$-ично иррациональной, если для всех $1\leq j\leq d,$ точка $x_j$ иррациональна относительна последовательности $P^j.$ \\

Для $\mathcal P$-ично иррациональных точек $x=(x_1,\ldots, x_d)\in R^d$ определим производную $\Psi'(x)$ и мажоранту $\Psi^*(x)$  следующим образом:
$$\Psi'(x)=\lim_{k\rightarrow\infty \atop{x\in I_k\in\Lambda^d_k}}\frac{\Psi(I_k)}{\mu(I_k)}, \ \Psi^*(x)=\sup_{I:x\in I\in\Lambda^d}\frac{|\Psi(I)|}{\mu(I)}.$$
\\\\

В настоящей работе нас будут интересовать $A\mathcal H$-интегралы по системам функций удовлетворяющих следующим условиям.\\
Зафиксируем некоторый $I_0\in\Lambda^d. $
 Пусть функции $h_i(x):I_0\rightarrow R$ такие, что
\begin{equation}\label{h1}
0\leq h_1(x)\leq h_2(x)\leq\dots\leq h_m(x)\leq\dots, \ \lim_{m\rightarrow\infty}h_m(x)=\infty
\end{equation} и существует постоянная $C>0$ и для каждого $m\in\mathbb N$ такие параллелепипеды $I^m_1,\dots,I^m_{n_m}\in\Lambda^d,$ что $\hbox{int}\,( I^m_i)\cap\hbox{int}\,( I^m_j)=\emptyset, i\neq j$
и $\cup_{k=1}^{n_m}I^m_k=I_0,$ для которых
\begin{equation}\label{h2}
\sup_{x\in I^m_k} h_m(x)\leq C\inf_{x\in I^m_k}h_m(x),\
\end{equation}
 для  всех $ m\in \mathbb N,\ 1\leq k\leq n_m,$ и
\begin{equation}\label{h3}
\inf_{m,k}\int_{I^m_k}h_m(x)dx>0.
\end{equation}
Иначе говоря, для каждой функции можно параллелепипед $I_0$ раздробить на маленькие "кусочки" \  , на каждом из которых эта функция принимает значения эквивалентные друг другу, a интегралы по этим "кусочкам"  \ больше некоторой положительной постоянной. \\

\begin{theo}\label{theorem psi}Пусть последовательность функций $h_n(x)$ удовлетворяет условиям \rf{h1},\rf{h2},\rf{h3}, а $\Psi$ комлекснозначная адитивная функция
определенная на $\mathcal P$-ичных параллелепипедах, где $\mathcal P=\{P^j\}_{j=1}^d, P^j=\{p^j_i\}_{i=1}^\infty$ и
\begin{equation}\label{M}
\sup_{i\in N, 1\leq j \leq d}p^j_i<\infty.
\end{equation}
 Если
\begin{equation}\label{int h}
\lim_{m\rightarrow\infty}\int_{\{x\in I_0;\Psi^*(x)>h_m(x)\}}h_m(x)dx=0
\end{equation}
и  п.в. существует $\Psi'(x),$ тогда для всякого $I\in\Lambda^d,$ $I\subset I_0,$ имеет место
$$\Psi(I)=\lim_{m\rightarrow\infty}\int_{I}\left[\Psi'(x)\right]_{h_m(x)}dx.$$
\end{theo}

\begin{proof}
Без ограничения общности можно считать, что $I=I_0=[0,1]^n.$
Пусть $I_{ki},$$\ i=1,2, \ldots,\prod_{j=1}^dm_k^j, $ суть все параллелепипеды ранга $k$ из $[0,1]^n.$
Обозначим
\begin{equation}\label{psi k def}
\Psi_k(x)=\frac{\Psi(I_{ki})}{\mu(I_{ki})}, \ \textrm{при} \ x\in \hbox{int}\,(I_{ki}).
\end{equation}
Ясно, что п.в. существует предел
\begin{equation}\label{psi k a.e.}
\lim_{k\rightarrow\infty}\Psi_k(x)=\Psi'(x)=:f(x).
\end{equation}
Введем следующие обозначения:
$$M=\sup_{i\in N, 1\leq j \leq d}p^j_i, \ \la^m_k=\inf_{x\in I^m_k}h_m(x)\ \textrm{ и} \ \ep_0=\inf_{m,k}\la^m_k\mu(I^m_k).$$
Из \rf{h2}, \rf{h3} вытекает, что $\ep_0>0.$
Возьмем $\ep<M^{-d}\ep_0.$
Из условия \rf{int h} следует, что для достаточно большого $m$ имеет место
\begin{equation}\label{int hm ep}
\int_{E_m}h_{m}(x)dx<\ep,\ \textrm{где}\ E_m=\{x\in[0,1]^d;\Psi^*(x)>h_{m}(x)\}.
\end{equation}
Следовательно из \rf{h2} имеем
\begin{equation}\label{laep}
\sum_{k=1}^{n_m}\la^m_k\mu\{x\in I^m_k,\Psi^*(x)>C\la^m_k\}<\ep.
\end{equation}
Зафиксируем $m$ и через $k_0$ обозначим наибольший ранг параллелепипедов $I^m_k,$ т.е.
$$k_0=\max_{1\leq k\leq n_m} r(I^m_k)
.$$
Обозначим $\tilde\Psi_{k_0}(x)=\Psi_{r(I^m_k)}(x),$ при $x\in I^m_k.$ Нетрудно убедиться, что $|\tilde\Psi_{k_0}(x)|\leq C\lambda^m_k,$ когда $x\in I^m_k.$
Действительно, в противном случае существует $k', 1\leq k' \leq n_m,$ такой что $|\tilde\Psi_{k_0}(x)|> C\lambda^m_{k'},$ для всех $x\in I^m_{k'},$ так как $\tilde\Psi_{k_0}(x)$ является постоянной
на каждом $I^m_k.$ Следовательно, из $\rf{laep}$ и определения $\ep_0$ будем иметь $\ep>\la^m_{k'}\mu(I^m_{k'})\geq\ep_0,$ что противоречит выбору $\ep.$\\
Пусть $I^m_1$  представлено в виде объединения параллелепипедов $I_1, \ldots, I_s$ ранга $r(I^m_1)+1.$
Если $|\Psi_{r(I^m_1)+1}(x)|$ не превосходит $C\la^m_1$ на каждом $I_j, j=1,\ldots, s,$ то на параллелепипеде $I^m_1$ положим $\tilde\Psi_{k_0+1}(x)=\Psi_{r(I^m_1)+1}(x),$ и каждый из параллелепипедов $I_1, \ldots, I_s$ назовем параллелепипедом первого класса для $\tilde\Psi_{k_0+1}(x).$ В противном случае положим $\tilde\Psi_{k_0+1}(x)=\Psi_{r(I^m_1)}(x)$ на $I^m_1$, и назовем $I^m_1$ параллелепипедом второго класса для $\tilde\Psi_{k_0+1}(x).$ Аналогично определим класс каждого параллелепипеда  $I^m_2,\ldots , I^m_{n_m},$ а также определим  $\tilde\Psi_{k_0+1}(x)$ на  каждом из этих параллелепипедов, т.е. на $I^m_2,\ldots , I^m_{n_m}.$ \\

Допустим $\tilde\Psi_{k_0+l}(x)$ определено. Определим $\tilde\Psi_{k_0+l+1}(x)$ следующим образом.
Параллелепипеды второго класса для $\tilde\Psi_{k_0+l}(x)$, останутся  параллелепипедами второго класса также для $\tilde\Psi_{k_0+l+1}(x), $ и положим $\tilde\Psi_{k_0+l+1}(x)=\tilde\Psi_{k_0+l}(x)$
на этих параллелепипедах. А если $I$ является параллелепипедом первого класса для $\tilde\Psi_{k_0+l}(x)$, то действуем следующим образом. Допустим $I$
представленно в виде объединения параллелепипедов $I_1, \ldots, I_s$ ранга $r(I)+1.$ Без ограничения общности будем считать, что $I\subset I^m_1.$ Если на каждом $I_j, j=1,\ldots, s$ модуль от $\Psi_{r(I)+1}(x)$ не превосходит $C\la^m_1,$ то положим $\tilde\Psi_{k_0+l+1}(x)=\Psi_{r(I)+1}(x)$ на $I,$ и каждый из параллелепипедов $I_1, \ldots, I_s$ назовем параллелепипедом первого класса для
$\tilde\Psi_{k_0+l+1}(x).$ В противном случае назовем $I$ параллелепипедом второго класса для
$\tilde\Psi_{k_0+l+1}(x)$, и положим  $\tilde\Psi_{k_0+l+1}(x)=\Psi_{r(I)}(x)$ на $I.$

Итак, функция $\tilde\Psi_{k_0+l}(x)$ является константой на $\mathcal P$-ичных параллелепипедах $I_1, \ldots, I_t$(вообще говоря, разных рангов) и
\begin{equation}\label{tildePsi}
\tilde\Psi_{k_0+l}(x)=\frac{\Psi(I_j)}{\mu(I_j)},\ \textrm{при}\ x\in\hbox{int}\,(I_j).
\end{equation}
По определению функции $\tilde\Psi_{k_0+l}(x)$ удовлетворяют
\begin{equation}\label{psi la}
|\tilde\Psi_{k_0+l}(x)|\leq C\lambda^m_k,\ \textrm{при}\ x\in\hbox{int}\,(I^m_k).
\end{equation}
Через $A_{k_0,l}$ обозначим объединение всех параллелепипедов, которые являются параллелепипедами второго класса для  $\tilde\Psi_{k_0+l}(x),$ и пусть $A_{k_0,l}^k=A_{k_0,l}\cap I^m_k.$ Докажем, что
\begin{equation}\label{Ck0kl}
\mu(A_{k_0,l}^k)\leq M^d\mu\{x\in I^m_k; \Psi^*(x)>C\la^m_k\}.
\end{equation}
Заметим, что множество $A_{k_0,l}^k$ это объединение параллелепипедов второго класса для  $\tilde\Psi_{k_0+l}(x),$ которые являются подмножеством  $I^m_k.$ Следовательно в каждом из этих параллелепипедов $I$ есть хоть один параллелепипед $J_I, $ такой что $r(J_I)=r(I)+1$ и $|\Psi_{r(J_I)}(x)|>C\la^m_k,$ при $x\in J_I.$ Следовательно $\Psi^*(x)>C\la^m_k,$ при $x\in J_I.$
Заметим также, что из определения величины $M$ следует, что $\mu(I)\leq M^d\mu(J_I),$ откуда  получим
$$\mu(A_{k_0,l}^k)=\sum\mu(I)\leq M^d\sum\mu(J_I)\leq M^d\mu\{x\in I^m_k; \Psi^*(x)>C\la^m_k\}.$$
Перейдем к оцениванию следующей величины:
$$\left|\Psi(I_0)-\int_{I_0}\left[\Psi'(x)\right]_{h_{m}(x)}dx\right|\leq
\sum_{k=1}^{n_m}\left|\int_{I^m_k}\tilde\Psi_{k_0+l}(x)dx-\int_{I^m_k\backslash E_m}\Psi'(x)dx\right|
+\int_{E_m}h_m(x)dx
$$
$$\leq \sum_{k=1}^{n_m}\left|\int_{(I^m_k\backslash E_m)\cap A_{k_0,l}^c}\left(\tilde\Psi_{k_0+l}(x)-\Psi'(x)\right)dx\right|+\int_{E_m}\left|\tilde\Psi_{k_0+l}(x)
\right|dx+
$$
$$\sum_{k=1}^{n_m}\int_{(I^m_k\backslash E_m)\cap A_{k_0,l}}\left(\left|\tilde\Psi_{k_0+l}(x)\right|+\left|\Psi'(x)\right|\right)dx+\int_{E_m}h_m(x)dx =I_1+I_2+I_3+I_4.$$
Из \rf{psi la} и определения $\la^m_k$ имеем, что $I_2\leq CI_4,$ следовательно из \rf{int hm ep} получим
\begin{equation}\label{I2I4}
I_2+I_4<(C+1) \ep.
\end{equation}
При $x\in I^m_k\backslash E_m$ имеем
\begin{equation}\label{Clamk}
|\Psi'(x)|\leq\Psi^*(x)\leq h_m(x)\leq C\la^m_k.
 \end{equation}
 Следовательно из  \rf{psi la} и \rf{Ck0kl} получим
$$
I_3\leq\sum_{k=1}^{n_m}2C\la^m_k\mu(A_{k_0,l}^k)\leq 2CM^d\sum_{k=1}^{n_m}\la^m_k\mu\{x\in I^m_k; \Psi^*(x)>C\la^m_k\},
 $$
откуда применив $\rf{laep}$ получим следующую оценку для $I_3$

\begin{equation}\label{I3}
I_3\leq2CM^d\ep.
\end{equation}
Осталось оценить $I_1.$ Из \rf{psi k a.e.} следует, что существует $l_0,$ такое что для всякого
$l'>l_0$ имеет место
\begin{equation}\label{Bk0l}
\mu\{x\in I^m_k, \left|\Psi_{k_0+l'}(x)-\Psi'(x)\right|\geq \ep\}<\frac{\ep}{\la^m_kn_m}.
\end{equation}
Выберeм $l$ так, чтобы для любого $k$ выполнялось неравенство $r(I^m_k)+l>k_0+l_0,$ для всех $k=1,2,\ldots n_m.$ Обозначим
$$B_{k_0,l}^k=\{x\in I^m_k;\left|\Psi_{r(I^m_k)+l}(x)-\Psi'(x)\right|\geq \ep\}, \ \ B_{k_0,l}=\cup_{k=1}^{n_m}B_{k_0,l}^k.$$
Заметим, что из \rf{Bk0l} следует, что $\mu(B_{k_0,l}^k)<\ep/(\la^m_kn_m),$ откуда,
учитывая, что  $\tilde\Psi_{k_0+l}(x)=\Psi_{r(I^m_k)+l}(x),$ при $x\in I^m_k\cap A_{k_0,l}^c, $ и неравенства \rf{psi la}, \rf{Clamk}, получим
$$
I_1\leq \sum_{k=1}^{n_m}\int_{((I^m_k\backslash E_m)\cap A_{k_0,l}^c)\backslash B_{k_0,l}}\left|\tilde\Psi_{k_0+l}(x)-\Psi'(x)\right|dx+
$$$$
\sum_{k=1}^{n_m}\int_{((I^m_k\backslash E_m)\cap A_{k_0,l}^c)\cap B_{k_0,l}}\left(\left|\tilde\Psi_{k_0+l}(x)\right|+\left|\Psi'(x)\right|\right)dx
\leq \ep\mu(I_0)+\sum_{k=1}^{n_m}2C\la^m_k\frac{\ep}{\la^m_kn_m}.$$
Итак имеем $I_1\leq (2C+1)\ep.$ Объединяя эту оценку вместе с оценками  \rf{I2I4}, \rf{I3}, будем иметь, что для достаточно больших $m$
$$\left|\Psi(I_0)-\int_{I_0}\left[\Psi'(x)\right]_{h_{m}(x)}dx\right|<(3C+2+2CM^d)\ep.$$
Теорема \ref{theorem psi} доказана.

\end{proof}

\begin{rem}\label{rem extend}
Любую  комлекснозначную адитивную функцию $\Psi$,
определенную на $\mathcal P$-ичных параллелепипедах, можно продолжить следующим образом.
Пусть $$I=\left[\frac{n_1}{m^1_{k_1}},\frac{n_1+1}{m^1_{k_1}}\right]\times\ldots\times
\left[\frac{n_d}{m^d_{k_d}},\frac{n_d+1}{m^d_{k_d}}\right],$$
тогда существуют $I_1,\ldots, I_s\in\Lambda^d,$ что $\hbox{int}\,(I_i)\cap \hbox{int}\,(I_j)=\emptyset,$
при $i\neq j$, и $I=\cup_{i=1}^sI_i.$ Положим $\Psi(I)=\sum_{i=1}^s\Psi(I_i).$ Очевидно, что результат теоремы \ref{theorem psi} останется верным для таких $I.$
\end{rem}
\begin{rem}
При выполнении условия \rf{M} имеем, что последовательность функций $\Psi_k(x)$ определенная в \rf{psi k def} является регулярным мартингалом. Следовательно согласно теоремы Ю. Шоу (см. \cite[с.~242]{G66} и
\cite{Ch62}) последовательность $\Psi_k(x)$ будет сходиться п.в. на множестве, где либо $\limsup_{k\rightarrow\infty}\Psi_k(x)<+\infty, $ либо $\liminf_{k\rightarrow\infty}\Psi_k(x)>-\infty,$ тем более на множестве, где $\sup_k|\Psi_k(x)|<+\infty.$  Откуда будем иметь, что при выполнении условия \rf{int h} последовательность $\Psi_k(x)$ сходится п.в., т.е. условие существования п.в. $\Psi'(x)$ в теореме \ref{theorem psi} и п.в. сходимость  ряда в теореме \ref{theo multy haar} ниже не существенно.

\end{rem}

Заметим, что взяв в теореме \ref{theorem psi} $p_i^j=2,$ для всех $i\in N,1\leq j\leq d,$  а функции
$h_m(x)\equiv \lambda_m$ получим следующую теорему 
 доказанную Г. Геворкяном в \cite{GG95}:
\begin{theo}\label{col lambda}
Пусть $\Phi$ адитивная функция определенная на двоичных кубах. Если существует последовательность
$\lambda_m\uparrow\infty$, что
\begin{equation}\label{lamda m}
\lim_{m\rightarrow\infty}\lambda_m\mu\{x\in I_0; \Phi^*(x)>\lambda_m\}=0,
\end{equation}
и п.в. существует $\Phi'(x),$ то для всякого двоичного куба $I\subset I_0$ имеет место следующая формула
$$\Phi(I)=\lim_{m\rightarrow\infty}\int_I\left[\Phi'(x)\right]_{\lambda_m}dx.$$
\end{theo}

Приведем пример аддитивной функции определенной на $[0,1]$, для которой удовлетворяются условия теоремы \ref{theorem psi} для некоторой последовательности функций $\mathcal H,$ однако, мажоранта которой не удовлетворяет условию \rf{lamda m} теоремы \ref{col lambda}, более того, невозможно разделить отрезок $[0,1]$ на несколько отрезков, на каждом из которых удовлетворялось \rf{lamda m}.

Пусть $\{\chi_n(x)\}_{n=1}^\infty $ система Хаара (см. \cite{KS}, стр. 70).
Принято также двухиндексная нумерация этой системы по следующему правилу. Если $n=2^k+i,$ $i=1,2,\ldots,2^k,$ $k\geq 0,$ то пологается $\chi_k^{(i)}(x):=\chi_n(x)$ и $\chi_0^{(0)}(x):=\chi_1(x).$ В таких обозначениях для $\chi_k^{(i)}(x)$ имеем
$$\chi_k^{(i)}(x)=\left\{
\begin{array}{lll}
0, &\ \hbox{при} & x\not\in[\frac{i-1}{2^k},\frac i{2^k}]; \\
2^{\frac k 2}, &\ \hbox{при} & x\in(\frac{i-1}{2^k},\frac {2i-1}{2^k}); \\
-2^{\frac k 2}, &\ \hbox{при} & x\in(\frac {2i-1}{2^k},\frac{i}{2^k}).
\end{array}
\right.$$
Значения функции $\chi_k^{(i)}(x)$ в точках разрыва равняется полусумме односторонних пределов, однако, для нас не важно значения в этих точках, так как в совокупности эти точки образуют множество меры нуль.
\begin{example}
Обозначим
$$k_n=\frac{n(n-1)}{2} \ \textrm{ и}\ \alpha(n,i)=\left(1-\frac1{2^{i-1}}\right)2^{k_n+i}+1.$$ Рассмотрим следующий ряд
\begin{equation}\label{haar series}
f(x)=\sum_{n=1}^\infty\sum_{i=1}^n2^{\frac{k_n+i}2}\chi_{k_n+i}^{\alpha(n,i)}(x)=
\sum_{n=1}^\infty a_n\chi_n(x).
\end{equation}
Фиксируем $j\in N.$
Ясно, что
$$2S^*f(x)\geq\left|2^{\frac{k_n+i}2}\chi_{k_n+i}^{\alpha(n,i)}(x)\right|,\ \textrm{для любых}\ i=1,\ldots,n,\ n\in N,$$
где $S^*(x)$ мажоранта частных сумм ряда \rf{haar series}. Следовательно
\begin{equation}\label{S^*}
\mu\left\{x\in\left[1-\frac1 {2^j},1\right];S^*(x)>2^m\right\}\geq
\end{equation}
$$\mu\left\{x\in\left[1-\frac1 {2^j},1\right];\left|2^{\frac{k_n+i}2}\chi_{k_n+i}^{\alpha(n,i)}(x)\right|>2^{m+1}\right\},
$$
$\ \textrm{для любых}\ i=1,\ldots,n,\ n\in N.$
Заметим также, что
\begin{equation}\label{supset}
\left[1-\frac1{2^{n-1}},1-\frac1{2^{n}}\right)
\supset \hbox{supp}\ \chi_{k_n+n}^{\alpha(n,n)}
\supset \hbox{supp}\ \chi_{k_{n+1}+n}^{\alpha(n+1,n)}\supset
\end{equation}$$
\supset \hbox{supp}\ \chi_{k_{n+2}+n}^{\alpha(n+2,n)}\ldots
\supset \hbox{supp}\ \chi_{k_{n+l}+n}^{\alpha(n+l,n)}\ldots\ .
$$
Из определения $k_n$ следует, что любое натуральное число, в частности $m+2,$ может быть представлено в виде $k_n+i, $ где $1\leq i\leq n.$
Обозначив $i'=\max(i,j+1),$ из \rf{supset} получим
$$\hbox{supp}\ \chi_{k_n+i'}^{\alpha(n,i')}(x)\subset \left[1-\frac1 {2^j},1\right] \
\textrm{и}\ \|2^{\frac{k_n+i'}2}\chi_{k_n+i'}^{\alpha(n,i')}\|_\infty=2^{k_n+i'}\geq 2^{m+2}.$$
Из последнего и \rf{S^*} получим, что
$$
\mu\left\{x\in\left[1-\frac1 {2^j},1\right];S^*(x)>2^m\right\}
\geq\mu\left(\hbox{supp}\ \chi_{k_n+i'}^{\alpha(n,i')}\right)=
\frac1{2^{k_n+i'}}\geq\frac1{2^{k_n+i+j}}.
$$
Следовательно
$$2^m\mu\left\{x\in\left[1-\frac1 {2^j},1\right];S^*(x)>2^m\right\}
\geq\frac1{2^{j+2}}.$$
Поэтому для любого $\lambda>1$ будем иметь
$\lambda\mu\left\{x\in\left[1-\frac1 {2^j},1\right];S^*(x)>\lambda\right\}
\geq\frac1{2^{j+3}},$ что означает
$$\liminf_{\lambda\rightarrow\infty}\lambda\mu\left\{x\in\left[1-\frac1 {2^j},1\right];S^*(x)>\lambda\right\}
>0,$$
для всех $j\in \mathbb N$.\\
Из этого следует, что невозможно отрезок $[0,1]$ разделить на конечное число отрезков $I_k,\ k=1,\ldots, s,$ так чтобы для любого $k$ имело место
$$\liminf_{\lambda\rightarrow\infty}\lambda\mu\left\{x\in I_k;S^*(x)>\lambda\right\}
=0,$$
так как можно выбрать $j$ настолько большим, что $\left[1-\frac1 {2^j},1\right]$ являлся подмножеством того отрезка $I_k$, которому принадлежит 1.
\\
Для всякого двоичного отрезка $I,$ т.е. отрезка вида $[m/2^k,(m+1)/2^k],$ обозначим
$$\Psi(I)=\lim_{N\rightarrow\infty}\int_I\sum_{n\leq N}a_n\chi_n(x)dx.$$

Правосторонний предел существует для всех двоичных интервалов $I,$ так как для достаточно большого
$N$(зависящего от $I$) интегралы от $n$-ой функции Хаара по $I$ будут равняться нулю, при всех $n\geq N$.  Заметим также, что для двоично иррациональной точки $x$
мы будем иметь, что $\Psi'(x)=f(x),$ а $\Psi^*(x)=S^*(x).$  Следовательно для каждого разбиения отрезка $[0,1]$ существует отрезок $I_k,$ что
$$
\liminf_{m\rightarrow\infty}\lambda\mu\{x\in I_k; \Psi^*(x)>\lambda\}>0.$$
\textsl{Поэтому теорема \ref{col lambda} не может быть применена для восстановления аддитивной функции $\Psi$  из $\Psi'(x)$.}\\

С другой стороны, назначив
$$h_n(x)=\left\{
\begin{array}{cc}
  2^{k_n+2}, & x\in\left[0,\frac12\right) \\
  2^{k_n+3}, & x\in\left[\frac12,\frac34\right) \\
  \vdots & \vdots \\
 2^{k_n+n+1}, & x\in\left[1-\frac1{2^{n-1}},1-\frac1{2^n}\right) \\
 2^n, & x\in \left[1-\frac1{2^{n}},1\right]
\end{array}
\right.$$
нетрудно заметить, что функции $h_n(x)$ удовлетворяют условиям \rf{h1}-\rf{h3}.
Так как  при $x\in\left[1-\frac1{2^{j-1}},1-\frac1{2^j}\right),\ 1\leq j\leq n$ имеет место
$$\sum_{n=1}^m\sum_{i=1}^n\left|2^{\frac{k_n+i}2}\chi_{k_n+i}^{\alpha(n,i)}(x)\right|\leq
2^{k_m+j+1}=h_m(x),$$
то использовав \rf{supset} получим, что
$$\{x\in[0,1];S^*(x)>h_m(x)\}\subset\bigcup_{i=1}^m\hbox{supp}\ \chi_{k_{m+1}+i}^{\alpha(m+1,i)}\cup \bigcup_{i=m+1}^\infty\hbox{supp}\ \chi_{k_{i}+i}^{\alpha(i,i)}.$$
Следовательно
$$\int_{\{x\in[0,1];S^*(x)>h_m(x)\}}h_m(x)dx\leq
 \sum_{i=1}^m\int_{\hbox{supp}\ \chi_{k_{m+1}+i}^{\alpha(m+1,i)}}h_m(x)dx
+\sum_{i=m+1}^\infty\int_{\hbox{supp}\ \chi_{k_{i}+i}^{\alpha(i,i)}}h_m(x)dx
$$
$$=\sum_{i=1}^m2^{k_m+i+1}\cdot\frac1{2^{k_{m+1}+i}}+2^m\sum_{i=m+1}^\infty\frac1{2^{k_{i+1}}}
\leq \frac{2m}{2^m}+\frac{2^{m+1}}{2^{k_{m+1}}}\rightarrow0.$$
\textsl{Последнее означает, что для  аддитивной функции $\Psi$ выполняются все условия теоремы \ref{theorem psi}, и следовательно она может быть восстановлена из $\Psi'(x)$ посредством $A\mathcal H$-интеграла.}
\end{example}

 Пусть $P=\{p_i\}_{i=1}^\infty$
последовательность натуральных чисел отличных от единицы, $m_0=1,\ m_j=\prod_{i=1}^jp_j.$ Тогда
$P$-ичное разложение точки $x\in[0,1)$  будет иметь вид
\begin{equation}\label{P}
x=\sum_{j=1}^\infty\frac {x_j}{m_j},\qquad 0\leq x_j\leq p_j-1,\ x_j\in Z,
\end{equation}
и для того, чтобы каждой точке $x\in[0,1)$ соответствовал единственный $P$-ичный ряд, для $P$-ично рациональных точек, т.е. для точек вида $l/m_n,$ где $l=0,1,\ldots, m_n-1$, выбираем суммы \rf{P}, которые состоят из конечного числа ненулевых слагаемых. Напомним определения системы   {\em  мультипликативных функций Прайса $\{\psi_k(x)\}_{k=0}^\infty$} и системы {\em обобщенных функций Хаара.}
Система Прайса определяется по формулам
$$\label{chi}\psi_k(x)=\exp\left(2\pi i\sum_{j=1}^n\frac{\alpha_jx_j}{p_j}\right),$$
где $k=\sum_{j=1}^n\alpha_jm_{j-1},\ 0\leq \alpha_j\leq p_j-1, \ j=1,2,\ldots, n.$
\\
А  система обобщенных функций Хаара $\{\chi_n^P(x)\}_{n=0}^\infty$ определяется по формулам

$$\chi_0(x)\equiv1,$$
а для $n\geq 1$
$$\chi_n(x)=\chi^k_{r,s}(x)=\left\{
\begin{array}{cc}
  \sqrt {m_k}e^{2\pi i x_{k+1}s/p_{k+1}}, &  x\in\left[\frac r{m_k};\frac {r+1}{m_k}\right) \\
  0, & x\in [0,1)\backslash\left[\frac r{m_k};\frac {r+1}{m_k}\right)
\end{array}
\right.,$$
где $n=m_k+r(p_{k+1}-1)+s-1,\ 0\leq r\leq m_k-1,\ 1\leq s\leq p_{k+1}-1.$\\

Заметим, что система Хаара соответствует последовательности $p_i=2,\ i\in N.$
Для всякого $\mathcal P=\{P_j\}_{j=1}^d$ рассмотрим ряд
\begin{equation}\label{haar}
\sum_{\overline{n}\in Z_+^d}a_{\overline{n}}\chi_{\overline{n}}^{\mathcal P}(x)=
\sum_{{n_j\in Z_+}\atop{j=1,\ldots, d}}a_{n_1\ldots n_d}\chi_{n_1}^1(x_1)\ldots\chi_{n_d}^d(x_d)
\end{equation}
и $\mathcal P$-ичные частные суммы этого ряда, т.е.
$$\sum_{\overline{n}<\widetilde{N}}a_{\overline{n}}\chi_{\overline{n}}^{\mathcal P}(x),$$
где под $\overline{n}<\overline{k}$ будем понимать $0\leq n_j<m^j_{k_j},$ для всех $j=1,2,\ldots,d,$ а $\widetilde{N}=(N,\ldots,N)\in R^d,$ для  любого натурального $N.$
Нетрудно заметить, что существует  ряд по системе Прайса
\begin{equation}\label{walsh}
\sum_{\overline{n}\in Z_+^d}b_{\overline{n}}\psi_{\overline{n}}^{\mathcal P}(x)
\end{equation}
 такой, что для всех $k\in N^d$ имеет место
\begin{equation}\label{overline n k}
\sum_{\overline{n}<\overline {k}}a_{\overline{n}}\chi_{\overline{n}}^{\mathcal P}(x)=
\sum_{\overline{n}< \overline {k}}b_{\overline{n}}\psi_{\overline{n}}^{\mathcal P}(x).
\end{equation}
В частности из этого следует, что $\mathcal P$-ичные частные суммы этого ряда совпадают $\mathcal P$-ичными частными суммами  ряда \rf{haar}:
\begin{equation}\label{S_N}
S_N^\chi(x):=\sum_{\overline{n}<\widetilde{N}}a_{\overline{n}}\chi_{\overline{n}}^{\mathcal P}(x)=
\sum_{\overline{n}<\widetilde{ N}}b_{\overline{n}}\psi_{\overline{n}}^{\mathcal P}(x)=:S_N^\psi(x).
\end{equation}
Итак, каждому ряду \rf{haar} можно сопоставить ряд \rf{walsh}, и наоборот.
Каждому ряду \rf{haar} поставим в соответствие аддитивную функцию $\Psi$ следующим образом

\begin{equation}\label{haar psi}
\Psi(I)=\lim_{N\rightarrow\infty}\int_I\sum_{\overline{n}<\widetilde{ N}}a_{\overline{n}}\chi_{\overline{n}}^{\mathcal P}(x)dx, \ \textrm{при}\ I\in\Lambda^d,
\end{equation}
Правосторонний предел существует для всех  $I\in\Lambda^d,$ так как для достаточно большого
$N$(зависящего от $I$) имеет место
$$\int_Ia_{\overline{n}}\chi_{\overline{n}}^{\mathcal P}(x)dx=0, \ \textrm{при} \ \max_{1\leq j\leq d}\left\{\frac{n_j}{m^j_N}\right\}\geq 1. $$

Из \rf{haar psi} будем иметь, что для каждой $\mathcal P$-ично иррациональной точки $x$ выполняются
\begin{equation}\label{psi' psi*}
\Psi'(x)=\lim_{N\rightarrow\infty}\sum_{\overline{n}<\widetilde{ N}}a_{\overline{n}}\chi_{\overline{n}}^{\mathcal P}(x) \ \textrm{и}\ \Psi^*(x)=\sup_N\left|\sum_{\overline{n}<\widetilde{ N}}a_{\overline{n}}\chi_{\overline{n}}^{\mathcal P}(x)\right|.
\end{equation}
Тогда имеет место следующая
\begin{theo}\label{theo multy haar}Пусть $\mathcal P$-ичные суммы кратного ряда Хаара \rf{haar} п.в. сходятся к функции
$f(x)$, a для $\mathcal P$ выполняется \rf{M}. Если существует последовательность функций $\{h_m(x)\}_{m=1}^\infty$ удовлетворяющая условиям \rf{h1}-\rf{h3}, для которой
$$
\lim_{m\rightarrow\infty}\int_{\{x\in [0,1]^d;\Psi^*(x)>h_m(x)\}}h_m(x)dx=0,$$
тогда для всякого $n=(n_1,\ldots,n_d)\in\mathbb{Z}_+^d$ имеем
$$a_{\overline{n}}=\lim_{m\rightarrow\infty}\int_{[0,1]^d}
\left[f(x)\overline{\chi_{\overline{n}}^{\mathcal P}(x)}\right]_{h_m^{\overline{n}}(x)}dx,$$
где $h_m^{\overline{n}}(x)=\|\chi_{\overline{n}}^{\mathcal P}(x)\|_\infty h_m(x).$

\end{theo}
{\bf Доказательство} аналогично рассуждениям для кратного ряда Хаара из \cite{GG95}.
Зафиксируем $\overline{n}\in N^d$. Функция $\chi_{\overline{n}}^{\mathcal P}(x)$ принимает постоянные значения на параллелепипедах $I_1,\ldots I_s$(не обязательно $\mathcal P$-ичных), равные по модулю
$\|\chi_{\overline{n}}^{\mathcal P}(x)\|_\infty$ и равняется нулю вне $\cup_{i=1}^sI_i.$ Для  $\mathcal P$-ичного параллелепипеда $I\subset I_i$ положим
\begin{equation}\label{Psi i}
\Psi_i(I)=\lim_{N\rightarrow\infty}\int_I\overline{\chi_{\overline{n}}^{\mathcal P}(x)}\sum_{\overline{k}<\widetilde N}a_{\overline{k}}\chi_{\overline{k}}^{\mathcal P}(x)dx, \quad i=1,2,\ldots, s.
\end{equation}
Ясно, что функции $\Psi_i(x)$ являются аддитивными на $\mathcal P$-ичных параллелепипедах $I_i,$
\begin{equation}\label{psi i '}
\Psi_i'(x)=\lim_{N\rightarrow\infty}\sum_{\overline{k}<\widetilde N}a_{\overline{k}}\chi_{\overline{k}}^{\mathcal P}(x)\overline{\chi_{\overline{n}}^{\mathcal P}(x)}=f(x)\overline{\chi_{\overline{n}}^{\mathcal P}(x)} \ \textrm{п.в.}
\end{equation}
и
\begin{equation}\label{psi i *}
\lim_{m\rightarrow\infty}\int_{\{x\in I_i;\Psi_i^*(x)>\|\chi_{\overline{n}}^{\mathcal P}(x)\|_\infty h_m(x)\}}\|\chi_{\overline{n}}^{\mathcal P}(x)\|_\infty h_m(x)dx=0.
\end{equation}
Применяя теорему \ref{theorem psi} и замечание \ref{rem extend}, из \rf{psi i '} и \rf{psi i *} получим
$$\Psi_i(I_i)=\lim_{m\rightarrow\infty}\int_{I_i}
\left[f(x)\overline{\chi_{\overline{n}}^{\mathcal P}(x)}\right]_{h_m^{\overline{n}}(x)}dx,\quad i=1,\ldots, s.$$
Суммируя по $i$ и приняв во внимание аддитивность интеграла и определение функций $\Psi_i$(см. \rf{Psi i}), мы получим
$$\lim_{m\rightarrow\infty}\int_{[0,1]^d}
\left[f(x)\overline{\chi_{\overline{n}}^{\mathcal P}(x)}
\right]_{h_m^{\overline{n}}(x)}dx=\sum_{i=1}^s\Psi_i(I_i)=
\lim_{N\rightarrow\infty}\int_{[0,1]^d}\overline{\chi_{\overline{n}}^{\mathcal P}(x)}
\sum_{\overline{k}<\widetilde N}a_{\overline{k}}\chi_{\overline{k}}^{\mathcal P}(x)dx$$

$$=\lim_{N\rightarrow\infty}\sum_{\overline{k}<\widetilde N}a_{\overline{k}}\int_{[0,1]^d}\chi_{\overline{k}}^{\mathcal P}(x)
\overline{\chi_{\overline{n}}^{\mathcal P}(x)}dx=a_{\overline{n}}.$$

Теорема доказана.\\

Заметим, что из равенства соответствующих $\mathcal P$-ичных частных сумм рядов \rf{S_N} и определения $\Psi(I)$ следует, что
$$\Psi(I)=\lim_{N\rightarrow\infty}\int_IS_N^\psi(x)dx,$$
а также для каждой $\mathcal P$-ично иррациональной точки $x$ из \rf{S_N} и \rf{psi' psi*}
будем иметь $\Psi'(x)=\lim_{N\rightarrow\infty}S_N^\psi(x)$ и  $\Psi^*(x)=\sup_N|S_N^\psi(x)|.$

Аналогичная теорема может  быть доказана также для системы Прайса. А именно имеет место следующая

\begin{theo}\label{theo multy walsh}Пусть $\mathcal P$-ичные частные суммы  ряда  \rf{walsh} п.в. сходятся к функции
$f(x)$, a для $\mathcal P$ выполняется \rf{M}. Если существует последовательность функций $\{h_m(x)\}_{m=1}^\infty$ удовлетворяющая условиям \rf{h1}-\rf{h3}, для которой
$$
\lim_{m\rightarrow\infty}\int_{\{x\in [0,1]^d;\Psi^*(x)>h_m(x)\}}h_m(x)dx=0,$$
тогда для всякого $n=(n_1,\ldots,n_d)\in\mathbb{Z}_+^d$ имеем
$$b_{\overline{n}}=\lim_{m\rightarrow\infty}\int_{[0,1]^d}
\left[f(x)\overline{\psi_{\overline{n}}^{\mathcal P}(x)}\right]_{h_m(x)}dx.$$

\end{theo}

{\bf Доказательство} аналогично рассуждениям для системы Прайса теоремы 1 из \cite{K03}.
Для каждого $P^j$ функции Прайса и соответствующие обобщенные функции Хаара одинакового ранга $n$ связаны друг с другом при помощи линейных зависимостей, т.е. для всякого   $k,\ m^j_{n-1}\leq k<m^j_n$ и $j,\ 1\leq j \leq d$ имеем
\begin{equation}\label{gamma}
\psi^j_k=\sum_{l=m_{n-1}}^{m_n-1}\gamma^{j,l}_k\chi_{l}^j.
\end{equation}
Покажем, что для каждого $j, 1\leq  j\leq d$ векторы $\Gamma_k^j=(\gamma_k^{j,m_{n-1}}, \ldots, \gamma_k^{j,m_n-1}),$ $m_{n-1}\leq k\leq m_n-1$ ортогональны друг другу воспользовавшись ортогональностью систем Прайса и обобщенного Хаара:
\begin{equation}\label{delta}
 \delta_{k,p}=(\psi_k^j,\overline{\psi_p^j})=(\sum_{l=m_{n-1}}^{m_n-1}\gamma^{j,l}_k\chi_{l}^j,
\sum_{s=m_{n-1}}^{m_n-1}\overline{\gamma^{j,s}_p\chi_{s}^j})
=\sum_{l,s=m_{n-1}}^{m_n-1}\gamma^{j,l}_k\overline{\gamma^{j,s}_p}(\chi_{l}^j,\overline{\chi_s^j})=
\end{equation}$$
\sum_{l,s=m_{n-1}}^{m_n-1}\gamma^{j,l}_k\overline{\gamma^{j,s}_p}\delta_{l,s}
=\sum_{l=m_{n-1}}^{m_n-1}\gamma^{j,l}_k\overline{\gamma^{j,l}_p.}
$$
Из \rf{overline n k} следует, что для любого $\overline s\in Z^d_+$
\begin{equation}\label{al bk}
\sum_{\overline s\leq \overline l<\overline s+\widetilde 1}a_{\overline l}\chi_{\overline l}^{\mathcal P}=\sum_{\overline s\leq \overline k<\overline s+\widetilde 1}b_{\overline k}\psi_{\overline k}^{\mathcal P},
\end{equation}
 а из  \rf{gamma} имеем
$$\sum_{\overline s\leq \overline k<\overline s+\widetilde 1}b_{\overline k}\psi_{\overline k}^{\mathcal P}=\sum_{\overline s\leq \overline k<\overline s+\widetilde 1}b_{k_1,\ldots,k_d}
\left(\sum_{l=m_{s_1}}^{m_{s_1+1}-1}\gamma_{k_1}^{1,l}\chi^1_l\right)\cdot\ldots\cdot
\left(\sum_{l=m_{s_d}}^{m_{s_d+1}-1}\gamma_{k_1}^{d,l}\chi^d_l\right)=$$
$$\sum_{\overline s\leq \overline k<\overline s+\widetilde 1}b_{\overline k}\cdot\left(\sum_{\overline s\leq \overline l<\overline s+\widetilde 1}\gamma^{\overline l}_{\overline k}\chi_{\overline l}^{\mathcal P}\right)=
\sum_{\overline s\leq \overline l<\overline s+\widetilde 1}\left(\sum_{\overline s\leq \overline k<\overline s+\widetilde 1}b_{\overline k}\gamma^{\overline l}_{\overline k}\right)\chi_{\overline l}^{\mathcal P},$$
где $\gamma^{\overline l}_{\overline k}=\gamma_{k_1}^{1,l_1}\cdot\ldots\cdot\gamma_{k_d}^{d,l_d}.$
Используя \rf{al bk}, получим следующее соотношение между коэффициентами $a_{\overline l}$ и $\ b_{\overline k}:$
$$a_{\overline l}=\sum_{\overline s\leq \overline k<\overline s+\widetilde 1}b_{\overline k}\gamma^{\overline l}_{\overline k}.$$
Отсюда, использовав равенства \rf{gamma}, \rf{delta} и теорему \ref{theo multy haar}, получим
$$b_{\overline p}=
\sum_{\overline s\leq \overline k<\overline s+\widetilde 1}b_{k_1,\ldots,k_d}\cdot \delta_{k_1,p_1}\cdot\ \ldots\ \cdot\delta_{k_d,p_d}
=\sum_{\overline s\leq \overline k<\overline s+\widetilde 1}b_{\overline k}\sum_{\overline s\leq \overline l<\overline s+\widetilde 1}\gamma^{\overline l}_{\overline k}\overline{\gamma^{\overline l}_{\overline p}}=$$
$$\sum_{\overline s\leq \overline l<\overline s+\widetilde 1} \overline{\gamma^{\overline l}_{\overline p}}\sum_{\overline s\leq \overline k<\overline s+\widetilde 1}b_{\overline k}\gamma^{\overline l}_{\overline k}=
\sum_{\overline s\leq \overline l<\overline s+\widetilde 1} \overline{\gamma^{\overline l}_{\overline p}}a_{\overline l}=
\sum_{\overline s\leq \overline l<\overline s+\widetilde 1} \overline{\gamma^{\overline l}_{\overline p}}\lim_{m\rightarrow\infty}\int_{[0,1]^d}
\left[f(x)\overline{\chi_{\overline{l}}^{\mathcal P}(x)}\right]_{h_m^{\overline{n}}(x)}dx=
$$
$$\sum_{\overline s\leq \overline l<\overline s+\widetilde 1} \overline{\gamma^{\overline l}_{\overline p}}\lim_{m\rightarrow\infty}\int_{[0,1]^d}
\left[f(x)\right]_{h_m(x)}\overline{\chi_{\overline{l}}^{\mathcal P}(x)}dx=
\lim_{m\rightarrow\infty}\int_{[0,1]^d}
\left[f(x)\right]_{h_m(x)}\sum_{\overline s\leq \overline l<\overline s+\widetilde 1}\overline{\gamma^{\overline l}_{\overline p}\chi_{\overline{l}}^{\mathcal P}(x)}dx=$$
$$=\lim_{m\rightarrow\infty}\int_{[0,1]^d}\left[f(x)\right]_{h_m(x)}\psi_{\overline{p}}^{\mathcal P}(x)dx=\lim_{m\rightarrow\infty}\int_{[0,1]^d}\left[f(x)\psi_{\overline{p}}^{\mathcal P}(x)\right]_{h_m(x)}dx.$$
Теорема доказана.\\

В последней теореме взяв $d=1$ и $h_m(x)\equiv \lambda_m$ получим теорему \ref{thK03} доказаную В. Костиным в \cite{K03}.

Автор благодарен  Г.Г. Геворкяну за обсуждение полученных результатов.
%
%
 %
%

\end{document}